\documentclass[12pt,sort&compress,preprint]{elsarticle}

\usepackage[english]{babel}
\usepackage[T1]{fontenc}
\usepackage[utf8]{inputenc}
\usepackage{lmodern}

\usepackage{amsmath,amsfonts,amssymb,amsthm}
\usepackage{mathtools}
\usepackage{commath}
\usepackage{esvect}

\usepackage{todonotes}

\usepackage{hyperref}

\newcommand{\ci}{\mathrm{i}}
\newcommand{\qi}{\mathbf{i}}
\newcommand{\qj}{\mathbf{j}}
\newcommand{\qk}{\mathbf{k}}
\newcommand{\eps}{\varepsilon}
\newcommand{\C}{\mathbb{C}}
\newcommand{\D}{\mathbb{D}}
\renewcommand{\H}{\mathbb{H}}
\newcommand{\R}{\mathbb{R}}
\renewcommand{\DH}{\D\H}
\newcommand{\cj}[1]{{#1}^\ast}
\newcommand{\ej}[1]{{#1}_\eps}

\newcommand{\Norm}[1]{N(#1)}
\newcommand{\SE}[1][3]{\operatorname{SE}(#1)}
\newcommand{\SO}[1][3]{\operatorname{SO}(#1)}

\DeclareMathOperator{\kernel}{kern}
\DeclareMathOperator{\rgcd}{realgcd}

\newtheorem{theorem}{Theorem}
\newtheorem{proposition}{Proposition}

\newtheorem{lemma}{Lemma}
\theoremstyle{definition}
\newtheorem{definition}{Definition}
\theoremstyle{remark}
\newtheorem{remark}{Remark}
\newtheorem*{remark*}{Remark}
\newtheorem{example}{Example}

\begin{document}

\begin{frontmatter}
\journal{} 

\title{Rational Motions of Minimal Quaternionic Degree with Prescribed Plane Trajectories}

\author{Zülal Derin Yaqub\corref{cor}}
\ead{zuelal.derin-yaqub@uibk.ac.at}
\author{Hans-Peter Schröcker}
\ead{hans-peter.schroecker@uibk.ac.at}
\address{University of Innsbruck, Department of Basic Sciences in Engineering Sciences, Technikerstraße~13, 6020 Innsbruck, Austria}

\cortext[cor]{Corresponding author}

\begin{abstract}
  This paper investigates the construction of rational motions of a minimal quaternionic degree that generate a prescribed plane trajectory (a “rational torse”). Using the algebraic framework of dual quaternions, we formulate the problem as a system of polynomial equations. We derive necessary and sufficient conditions for the existence of such motions, establish a method to compute solutions and characterize solutions of minimal degree. Our findings reveal that a rational torse is realizable as a trajectory of a rational motion if and only if its Gauss map is rational. Furthermore, we demonstrate that the minimal degree of a motion polynomial is geometrically related to a drop of degree of the Gauss and algebraically determined by the structure of the torse's associated plane polynomial and the real greatest common divisor of its vector part. The developed theoretical framework has potential applications in robotics, computer-aided design, and computational kinematic, offering a systematic approach to constructing rational motions of small algebraic complexity.
\end{abstract}

\begin{keyword}
  kinematics \sep dual quaternion \sep rational kinematic torse \sep developable
  surface \sep quaternionic polynomial
  \MSC[2020]{70B10, 51J15, 51N15, 70E15, 14J26, 11R52, 65D17  }
\end{keyword}

\end{frontmatter}

\section{Introduction}
\label{sec:introduction}

It is common to design a rigid body motion by prescribing one of its point
trajectories \cite{artobolevskii64,Wagner97,proskova17}. But one can also
imagine situations when the trajectory of a plane is to be created by a motion.
A single point trajectory curve does not constrain at all the orientation along
the curve while a plane trajectory partially constrains position and
orientation. Thus, there is an infinity of motions along a certain point or
plane trajectory and one can ask the question of how to sensible use these
degrees of freedom?

In this paper we consider rational rigid body motions of Euclidean three-space
with a prescribed rational plane trajectory (a \emph{rational torse} or
\emph{rational developable surface}).\footnote{We prefer the term ``torse'' over
  ``developable surface'' as it emphasizes that we consider one-parametric sets
  of planes (curves in dual space), not their envelopes (developable surfaces).}
We provide answers to the following questions:
\begin{itemize}
\item Under which conditions on the rational torse do such rigid body motions
  exist?
\item What is the minimal possible motion degree that can be achieved?
\item Under which conditions is this minimal motion essentially unique?
\end{itemize}
Here, the notion of \emph{degree} refers to the representation of rigid body
motions in the dual quaternion model of space kinematics (the ``quaternion
degree'' in the sense of \cite{juettler93}). It is different from the more
common degree with respect to the parametrization of $\SE$ by homogeneous
transformation matrices \cite{juettler93}. However, our notion of degree is
arguably more relevant when it comes to the construction of mechanisms to
perform a prescribed motion \cite{Gallet16,li18}. ``Essential uniqueness'' means
uniqueness in a geometric sense, up to coordinate changes in the moving frame,
that is, up to right multiplication with constant rigid body displacements.

While point trajectories in kinematics have been investigated extensively, much
less in known for the trajectory of planes (and of lines). A standard reference
on theoretical kinematics features a short section on that topic
\cite[Section~VI.7]{Bottema90}. The paper \cite{Selig00} studies rigid body
transformations of points, lines, and planes in a Clifford algebra setting and
discusses examples from computer vision and geometric constraint solving. An
additional motivation to study the kinematics of plane trajectories is the close
relation between torses and developable surfaces \cite{Bodduluri93}.

Our research questions are motivated by a similar investigation in \cite{li16} on minimal degree motions with a prescribed point trajectory $\gamma(t)$. The main results of that paper are as follows:
\begin{itemize}
\item The minimal degree motion to a rational curve is essentially unique in all
  cases.
\item In generic cases the minimal achievable motion degree equals the curve degree and the minimal degree motion is the trivial curvilinear translation along the curve.
\item If the curve is of circularity $c > 0$ (refer to \cite{li16} for a definition of that concept) the minimal motion degree is $\deg\gamma(t)-c$ and the motion itself is the superposition of a translation along the curve $\gamma(t)$ with some change of orientation.
\end{itemize}

The minimal degree motion was used in \cite{li18} for a construction of
exceptionally simple Kempe linkages for rational curves. In similar spirit, we
may ask for linkages to produce a rational plane trajectory. This is an open
question and its solution will require properties of rational motions of minimal
degree with a given plane trajectory. Results of this article will show that
there are substantial differences to the curve case:
\begin{itemize}
\item Rational motions with prescribed rational torse exist if and only if the
  torse has a rational Gauss image $n$. This is not always the case.
\item Generically, the minimal achievable motion degree is half the torse
  degree. It increases if there is a drop in the degree of the Gauss map~$n$.
\item The motion of minimal degree is essentially unique if and only if the
  degree of the Gauss map does not drop.
\end{itemize}

In case of prescribed point trajectories the essential quantities to determine the minimal motion degree are the curve's degree and its circularity while in case of prescribed plane trajectories it is the torse's degree and the degree of its Gauss image.

The remainder of this article is structured as follows:
Section~\ref{sec:preliminaries} provides the necessary mathematical background
on the dual quaternion model of space kinematics as well as on rational torses
and rational motions. It also formalizes the problem statement and reduces it to
the discussion of solutions of a system of equations over the ring of quaternion
polynomials. In Section~\ref{sec:system} we extensively discuss this system
of equations and provide results on its solubility, uniqueness of solutions and
solution degree. Our proof in that section are mostly constructive and can be
used to actually compute rational motion with a given rational plane trajectory.
Our main results are formulated and proved in Section~\ref{sec:minimal}.

\section{Preliminaries}
\label{sec:preliminaries}

One of the most important applications of the skew field $\H$ of quaternions is a
parametrization of the special orthogonal group $\SO$. Here, we describe this
well-known construction, we extend it to an equally well-known parametrization
of $\SE$ via dual quaternions, and we describe how to use it for the
transformation of planes. An extension of scalars from real numbers $\R$ to
rational functions $\R(t)$ then leads to a parametric version of this action as
well as to the description of rational motions via motion polynomials and of
rational torses via plane polynomials.

\subsection{Quaternions and Spherical Kinematics}

As usual, we write a quaternion $q \in \H$ as $q = q_0 + q_1\qi + q_2\qj +
q_3\qk$ with $q_0$, $q_1$, $q_2$, $q_3 \in \R$. The real number $q_0$ is called
the quaternion's \emph{scalar part} while $\vec{q} \coloneq q_1\qi + q_2\qj +
q_3\qk$ is its \emph{vector part.} Both scalar and vector part can be written in
terms of the \emph{conjugate quaternion} $\cj{q} \coloneqq q_0 - q_1\qi - q_2\qj
- q_3\qk$:
\begin{equation*}
  q_0 = \tfrac{1}{2}(q + \cj{q}),\quad
  \vec{q} = \tfrac{1}{2}(q - \cj{q}).
\end{equation*}
Quaternion conjugation satisfies the important rule $\cj{(ab)} = \cj{b}\cj{a}$
for any $a$, $b \in \H$. A quaternion with zero scalar part is called
\emph{pure} or \emph{vectorial.} The quaternion \emph{norm} is the non-negative
real number $\Norm{q} \coloneqq q\cj{q} = q_0^2 + q_1^2 + q_2^2 +
q_3^2$.\footnote{The name ``norm'' is common for quaternion algebras but
  diverges from the usual concept of a ``norm'' in linear algebra where
  $\Norm{q}$ would rather be referred to as ``squared norm.''} It is
multiplicative, ie. $N(ab) = N(a)N(b)$ and it gives rise to the
inverse quaternion
\begin{equation*}
  q^{-1} = \frac{\cj{q}}{\Norm{q}},
\end{equation*}
provided $q \neq 0$.

Identifying the point $x = (x_1,x_2,x_3) \in \R^3$ with the vectorial quaternion $x_1\qi + x_2\qj + x_3\qk$, we define an action of $q \in \H \setminus \{0\}$ on $\R^3$ as
\begin{equation}
  \label{eq:1}
  x \mapsto qxq^{-1} = \frac{q x \cj{q}}{\Norm{q}}.
\end{equation}
It is the rotation with axis $\vec{q}$ and rotation angle $\varphi$ given by
\begin{equation*}
  \tan\frac{\varphi}{2} =
  \frac{\sqrt{\vec{q}\cj{\vec{q}\,}}}{q_0}.
\end{equation*}
Importantly, composition of rotations of the shape \eqref{eq:1} corresponds to
multiplication of quaternions. In other words, the map that sends the quaternion
$q$ to the rotation \eqref{eq:1} is a homomorphism between the multiplicative
group of quaternions and $\SO$. Since precisely the non-zero real multiples of
$q$ describe the same rotation, this homomorphism provides an \emph{isomorphism}
between $\SE$ and the factor group $\H^\times/\R^\times$, where $\H^\times$ and
$\R^\times$ denote the quaternionic and the real multiplicative groups,
respectively.

\subsection{Dual Quaternions and Space Kinematics}

With the aim of also incorporating translations into this multiplicative action, we consider dual quaternions $\DH$. They are obtained from quaternions by extensions of scalars from real numbers $\R$ to dual numbers $\D \coloneqq \R[\eps]/{\langle \eps^2 \rangle}$. A dual number is given as $a + \eps b$ where $a$ and $b$ are real; multiplication is determined by the rule that $\eps$ squares to zero, i.e., $(a+b\eps)(c+d\eps) = ac + (ad+bc)\eps$.

A dual quaternion $h \in \DH$ can be written as $h = p + \eps d$ where \emph{primal part} $p$ and \emph{dual part} $d$ are both quaternions. We will use its vector part $\vec{d} = \vec{p} + \eps \vec{d}$, its conjugate $\cj{h} = \cj{p} + \eps\cj{d}$, and its $\eps$-conjugate $\ej{h} \coloneqq p - \eps d$. The dual quaternion norm
\begin{equation*}
  \Norm{h} \coloneqq h\cj{h} = p\cj{p} + \eps(p\cj{d} + d\cj{p})
\end{equation*}
is no longer a real but a dual number as both, $p\cj{p}$ and $p\cj{d} +
d\cj{p}$, are real. With $p = p_0 + p_1\qi + p_1\qj + p_1\qk$ and $d = d_0 +
d_1\qi + d_1\qj + d_1\qk$ the vanishing of the dual part can be expressed
explicitly as
\begin{equation}
  \label{eq:2}
  p\cj{d} + d\cj{p} = p_0d_0 + p_1d_1 + p_2d_2 + p_3d_3 = 0.
\end{equation}
Equation~\eqref{eq:2} is commonly referred to as the \emph{Study condition.} Now
denote by $\DH^\times$ the multiplicative group
\begin{equation*}
  \DH^\times \coloneqq \{ h = p + \eps d \in \DH \mid h\cj{h} \in \R^\times \}.
\end{equation*}
It consists of dual quaternions $h$ of non-zero real norm whose inverse is given
by
\begin{equation*}
  h^{-1} = \frac{\cj{h}}{\Norm{h}}.
\end{equation*}

Similar to the isomorphism $\H^\times/\R^\times \to \SO$, there is an
isomorphism $\DH^\times/\R^\times \to \SE$ but its construction requires a
slightly different embedding of points. Identify $x = (x_1,x_2,x_3)$ with the
dual quaternion $1 + \eps(x_1\qi + x_2\qj + x_3\qk) = 1 + \eps\vec{x}$ and
define the action of $h = p + \eps d \in \DH^\times$ as
\begin{equation}
  \label{eq:3}
  x \mapsto \ej{h} x h^{-1}
  = \frac{\ej{h} x \cj{h}}{\Norm{h}}
  = \frac{p\cj{d} - d\cj{p} + \eps p\vec{x}\cj{p}}{\Norm{p}}.
\end{equation}
The shape of Equation~\eqref{eq:3} clearly shows that it is the composition of a
rotation, described by $p$, and a translation, encoded by both $p$ and~$d$.

In this article we prefer a homogeneous formulation of the thus provided
isomorphism between $\DH^\times / \R^\times$ and $\SE$. Describing points by
homogeneous coordinate vectors $[x_0,x_1,x_2,x_3]$ and identifying them with $x
= x_0 + \eps(x_1\qi + x_2\qj + x_3\qk) \in \DH$, the action \eqref{eq:3} simply
becomes
\begin{equation}
  \label{eq:4}
  x \mapsto \ej{h}x\cj{h}.
\end{equation}
It is preferable in the context of rational kinematics as it avoids the
square roots when normalizing~$h$.

The action \eqref{eq:4} on points implicitly defines an action on planes. In
order to describe it explicitly, we embed the plane with equation $u_0 + u_1x+
u_2y + u_3z = 0$ into the dual quaternions as $u = u_1\qi + u_2\qj + u_3\qk +
\eps u_0 = \vec{u} + \eps u_0$. Then, the dual quaternion $h = p + \eps d \in
\DH^\times$ acts on this plane according to
\begin{equation}
  \label{eq:5}
  u = \vec{u} + \eps u_0 \mapsto \ej{h} u \cj{h} = p\vec{u}\cj{p} + \eps(p\cj{d}-d\cj{p})
\end{equation}
\cite{Selig00}. This action is well-defined as $p\vec{u}\cj{p}$ is indeed
vectorial and $p\cj{d}-d\cj{p}$ is a scalar.

\subsection{Rational Curves, Torses, and Motions}

Central concepts of this article are rational torses (rational curves in
projective dual space) and rational motions. We may define them by yet another
extension of scalars, considering dual quaternions not just over the real
numbers $\R$ but over the field $\R(t)$ of rational functions in one
indeterminate $t$ that will also attain the role of a parameter for the
parametric torse or motion. The indeterminate $t$ commutes with all dual
quaternions and is unaffected by conjugation.

Now, we define formally:

\begin{definition}
  \label{def:rational}
  \emph{Rational curves, rational torses,} and \emph{rational motions} are points in the projective space over the vector space $\DH$ of dual quaternions with base field $\R(t)$, subject to the following conditions:
  \begin{itemize}
  \item A rational curve $[\gamma]$ has vanishing primal vector part and dual scalar part, that is $\gamma = \gamma_0 + \eps(\gamma_1\qi + \gamma_2\qj + \gamma_3\qk)$.
  \item A rational torse ${[u]}$ has vanishing primal scalar part and dual vector part, that is $u = u_1\qi + u_2\qj + u_3\qk + \eps u_0$.
  \item For a rational motion ${[C]}$, the Study condition \eqref{eq:2} is satisfied and the primal part is not zero, that is, $C = P + \eps D$ where $P \neq 0$ and $P\cj{D} + D\cj{P} = 0$.
  \end{itemize}
\end{definition}

We use square brackets around $\gamma$, $u$, and $C$ to denote points of the
respective projective spaces that these entities lie in. It is clear from
Definition~\ref{def:rational} that rational curves, torses and motions can be
represented not just by rational functions but also by polynomials in the
indeterminate $t$ and with coefficients from $\DH$. This is what we will usually
do. If $[C]$ is a polynomial representation of a rational motion, we call $C$ a
\emph{motion polynomial}, cf. \cite{hegedus13}. If $[u]$ is a polynomial
representation of a rational torse, we call $u$ a \emph{plane polynomial.} The
extensions of the actions \eqref{eq:3} and \eqref{eq:5} from dual quaternions to
motion polynomials provides polynomial descriptions for rational curves and
torses, respectively. Note that we also constant points, plans, and rigid body
displacements are to be considered as points of the respective projective
spaces.

An important tool in some of our proofs is polynomial division in the algebra of
quaternion polynomials. Quaternionic polynomial division works similar to the
division of real or complex polynomials but, due to non-commutativity, comes in
a left and in a right version:

\begin{proposition}
  \label{prop:division}
  Given polynomials $F$, $G \in \DH[t]$ with the leading coefficient of $G$
  being invertible, there exist unique polynomials $Q_\ell$, $Q_r$, $R_\ell$,
  $R_r$ (left/right quotients and remainders) such that $\deg R_\ell < \deg G$,
  $\deg R_r < \deg G$, and $F = Q_\ell G + R_\ell = Q_r G + R_r$. If the divisor
  $G$ is real, left- and right-division produce the same quotients and
  remainders.
\end{proposition}

Proposition~\ref{prop:division} is well-known, also in the context of more
general rings. A proof for quaternionic polynomials is given in
\cite[Proposition~4]{Damiano10}. Its extension to dual quaternion polynomials is
straightforward and the only essential additional requirement is that the
divisor's leading coefficient is invertible \cite{li18b}. More information on
quaternionic polynomials in general can be found at diverse places, for example,
\cite{Falcao17}.

\subsection{Problem Statement and Simplifications}
\label{sec:problem-statement}

Given a rational torse $[u] = [\vec{u} + \eps u_0] = [u_1\qi + u_2\qj + u_3\qk +
\eps u_0]$, we are looking for a rational motion $[C] = [P + \eps D]$ and a
constant plane $[w]$ (the ``moving plane'') such that $[u] = [\ej{C} w \cj{C}]$
holds. For any plane $[w]$ there exists a suitable rigid body displacement $a$
such that $\ej{a}v\cj{a} = \qk$. The dual quaternion $a$ may be absorbed in
$[C]$ (replace $C$ by $Ca$) so that we can assume without loss of generality
that the moving plane is $[\qk]$. We thus will study the equation
\begin{equation}
  \label{eq:6}
  [u] = [\ej{C}\qk\cj{C}].
\end{equation}
The rational torse $[u]$ is given, the rational motion $[C]$ is sought. As we
shall see (Theorem~\ref{th:1}), solutions to \eqref{eq:6} exist if and only if a
natural assumption on $[u]$, rationality of its Gauss map, is satisfied. If a
solution exists, it is not unique: We can always multiply $[C]$ from the right
with a rational motion $[E]$ that fixes the plane $[\qk]$. The motion $[E]$ is
then a \emph{planar} motion, that is, a motion in $\SE[2]$, that necessarily is
of the shape $E = e_0 + e_3\qk + \eps(e_5\qi + e_6\qj)$. Indeed, $\ej{E} \qk
\cj{E} = (e_0^2 + e_3^2)\qk$ whence $[\ej{E} \qk \cj{E}] = [\qk]$. Given the
non-uniqueness of solutions to Equation~\eqref{eq:6}, it is natural to search
for ``simple'' solutions and it is equally natural to define simplicity in terms
of the motion degree:

\begin{definition}
  \label{def:degree}
  The \emph{degree} of a rational curve, torse, or motion is the minimal degree of a \emph{polynomial} that represents the respective curve, torse, or motion.
\end{definition}

Note that Definition~\ref{def:degree} refers to the dual quaternion model we use. It does not really matter when talking about curves and torses but our notion of motion degree differs from the more common degree of rational motions that are given in terms of homogeneous transformation matrices, cf. \cite{juettler93,li16} for more information on that subject.

\begin{example}
  \label{ex:1}
  The rational torse
  \begin{equation*}
    [u] = \bigl[ \tfrac{t+1}{t}\qi + \tfrac{1}{t-1}\qj +
    \tfrac{2(t+1)}{t-1}\qk + \tfrac{t+1}{t(t-1)}\eps \bigr]
  \end{equation*}
  is of degree two because of
  \begin{multline*}
    [u] =
\bigl[ \tfrac{1}{t(t-1)}((t+1)(t-1)\qi + t\qj + 2t(t+1)\qk +
    (t+1)\eps) \bigr] \\=
    [ (t+1)(t-1)\qi + t\qj + 2t(t+1)\qk + (t+1)\eps ].
  \end{multline*}
\end{example}

Call a polynomial in $w = w_0 + w_1\qi + w_2\qj + w_3\qk + \eps(w_4 + w_5\qi +
w_6\qj + w_7\qk) \in \DH[t]$ \emph{reduced} if it has no real polynomial factor
of positive degree. This can be expressed using the $\gcd$ of real polynomials.
For $w \in \DH[t]$, define its real greatest common divisor as $\rgcd(w)
\coloneqq \gcd(w_0,w_1,\ldots,w_7)$. Now, $w$ is reduced if and only if
$\rgcd(w) = 1$. The degree of a rational curve, torse, or motion $[x]$ equals
$\deg x$ if and only if $x$ is reduced.

The motion polynomial $C \in \DH[t]$ solves \eqref{eq:6} if and only if there
exists a polynomial $h \in \R[t]$ such that
\begin{equation}
  \label{eq:7}
  \ej{C} \qk \cj{C} = h u.
\end{equation}
It is this polynomial version of \eqref{eq:6} that we will mostly study in the
remainder of this article, most notably in Section~\ref{sec:system}.

Since real polynomial factors of $u$ can be absorbed into $h$, it seems
reasonable to assume that $u = u_1\qi + u_2\qj + u_3\qk + \eps u_0$ is a reduced
polynomial. However, we will not generally do this and instead rely on the more
subtle notion of \emph{saturated} or \emph{minimally saturated} plane
polynomials (cf. Definition~\ref{def:saturated-polynomial} below). A general
technical assumption that we would like to make is
\begin{equation*}
  \deg u_0 = \deg u_1 = \deg u_2 = \deg u_3.
\end{equation*}
Should this not be fulfilled, we can apply a suitable re-parametrization of the
shape
\begin{equation*}
  t \mapsto \tfrac{\alpha t + \beta}{\gamma t + \delta}
\end{equation*}
with some $\alpha$, $\beta$, $\gamma$, $\delta \in \R$ such that $\alpha\delta -
\beta\gamma \neq 0$ and multiply away denominators. An alternative would be the
use of homogeneous polynomials throughout this text. We do not do this because
of the overhead in notation and also existing literature on motion polynomials
avoids this.

\section{A System of Equations for Quaternionic Polynomials}
\label{sec:system}

Given a plane polynomial $u = \vec{u} + \eps u_0 \in \DH[t]$ this section
presents a thorough discussion of solutions to Equation~\eqref{eq:7} for a real
polynomial $h \in \R[t]$ and a motion polynomial $C \in \DH[t]$. Typically we
assume that $h$ is given and satisfies some specific properties and we try to
show existence or non-existence of a suitable $C$. Equation~\eqref{eq:7} is
obviously related to the computation of rational motions with a prescribed plane
trajectory.

\subsection{The Primal Part}
\label{sec:primal}

To begin with, we study the primal part of Equation~\eqref{eq:7} in the special
case $h = 1$:
\begin{equation}
  \label{eq:8}
  P \qk \cj{P} = \vec{u}.
\end{equation}
Taking norms on both sides of \eqref{eq:8} gives
\begin{equation*}
  \Norm{P}^2 = \Norm{\vec{u}}.
\end{equation*}
Since $\vec{u}$ and $u$ have the same norm, it is necessary for solubility of
\eqref{eq:8} that the norm of $u$ is a square. This motivates the following
definition:
\begin{definition}
  \label{def:kinematic}
  The plane polynomial $u = \vec{u} + \eps u_0 \in \DH[t]$ is called
  \emph{kinematic} if its norm is a square in $\R[t]$.
\end{definition}

We continue with two lemmas that study the possibilities to create a real
polynomial $g$ by complex polynomial conjugation of the fixed quaternion $\qk$
and to create a vectorial polynomial $\vec{u}$ by quaternionic polynomial
conjugation of $\qk$. Both are important technical ingredients in subsequent
proofs and constructions but not really new. Consequently, a large part of our
``proofs'' will consist of pointers to literature.

In what follows, we embed the complex number field $\C$ into $\H$ as sub-algebra
generated by $1$ and~$\qk$.

\begin{lemma}
  \label{lem:1}
  Given a monic polynomial $g \in \R[t]$ there exists a polynomial $G \in \C[t]$
  such that $G\qk\cj{G} = \ell g \qk$ where $\ell \in \R[t]$. If
  \begin{equation*}
    g = \prod_{i=1}^m (t - t_i)^{\mu_i}
    \prod_{i=1}^n (t - a_i - b_i\qk)(t - a_i + b_i\qk)
  \end{equation*}
  with pairwise different real values $t_1$, $t_2$, \ldots, $t_m$, positive
  integers $\mu_i$, and not necessarily different complex numbers $a_1 +
  b_1\ci$, $a_2 + b_2\ci$, \ldots, $a_n + b_n\ci$ is the factorization of $g$
  over $\C$, then suitable polynomials $G$ and $\ell$ of minimal degrees are
  \begin{equation}
    \label{eq:9}
    G = \prod_{i=1}^m (t - t_i)^{\lceil \mu_i/2 \rceil} \prod_{i=1}^n (t - a_i - b_i\qk),
    \quad
    \ell = \prod_{i=1}^m (t - t_i)^{\mu_i - \lceil \mu_i/2 \rceil}
  \end{equation}
  The polynomial $\ell$ of minimal degree is unique, the polynomial $G$ is
  unique up to the conjugation of factors.
\end{lemma}

We omit the simple technical proof of Lemma~\ref{lem:1} but illustrate it at hand of an example:

\begin{example}
  \label{ex:2}
  Consider the polynomial $g = t^2(t-1)^3(t^2+1)$. The Equation~\eqref{eq:9}
  suggests $G = t(t-1)^2(t-\qk)$. Indeed,
  \begin{equation*}
    G\qk\cj{G} = t(t-1)^2(t-\qk)\qk t(t-1)^2(t+\qk)
    = t^2(t-1)^4(t^2+1)\qk
    = \ell g \qk.
  \end{equation*}
  where $\ell = t-1$. Alternatively, we could have used $G = t^2(t-1)(t + \qk)$
  and the same~$\ell$.
\end{example}

The gist of Lemma~\ref{lem:1} is that irreducible real quadratic factors and
linear real factors with \emph{even} multiplicity can appear in products of the
form $G\qk\cj{G}$ with $G \in \C[t]$.

\begin{lemma}
  \label{lem:2}
  Given a vectorial polynomial $\vec{u} \in \H[t]$, set $g\coloneqq
  \rgcd(\vec{u})$ and $\vec{v} \coloneqq \vec{u}/g$. By Lemma~\ref{lem:1} there
  exist $\ell \in \R[t]$ and $G \in \C[t]$, both of minimal, such that
  $G\qk\cj{G} = \ell g \qk$. If $\vec{v}$ is kinematic, there exists a
  polynomial $Q \in \H[t]$ such that with $P \coloneqq QG$ we have $P\qk\cj{P} =
  \ell\vec{u}$. The polynomial $Q$ is unique up to right multiplication with a
  complex number of unit norm and has no real polynomial factors.
\end{lemma}

\begin{proof}
  Because of $\rgcd(\vec{v}) = 1$, \cite[Theorem~4.2]{choi02} or
  \cite[Lemma~2.3]{schroecker24} imply existence of $Q \in \H[t]$ such that
  $\vec{v} = Q\qk\cj{Q}$. With $P = QG$ we then clearly have
  \begin{equation*}
    P\qk\cj{P} = (QG)\qk\cj{(QG)} = Q(G\qk\cj{G})\cj{Q} = \ell g Q\qk\cj{Q} = \ell g \vec{v} = \ell \vec{u}.
  \end{equation*}
  Uniqueness of $Q$ up to right multiplication with unit complex numbers is a
  special case of the uniqueness statement in \cite[Theorem~2]{li16}. Finally,
  real polynomial factors of $Q$ are not possible because they would also be
  factors of $\vec{v}$ which contradicts $\rgcd(\vec{v}) = 1$.
\end{proof}

\begin{remark*}
  The statement of Lemma~\ref{lem:2} is well-known in the context of curves with
  a Pythagorean hodograph \cite{farouki08}. The proofs of
  \cite[Theorem~4.2]{choi02} and of \cite[Lemma~2.3]{schroecker24} are both
  constructive and allow to actually compute the polynomial $Q \in \H[t]$. In
  this article we never actually perform this computation and just rely on the
  fact that a suitable quaternionic polynomial $Q$ exists.
\end{remark*}

Lemmas~\ref{lem:1} and \ref{lem:2} indicate that not all real polynomial factors
can appear on the right-hand side of \eqref{eq:8}. In fact, conjugation of $\qk$
with a polynomial $P \in \H$ can only produce vectorial quaternionic polynomials
whose real $\gcd$ factors over $\R$ into a product of quadratic irreducible
polynomials and linear polynomials of \emph{even} multiplicity. Hence, there
should be some restriction on the vector part $\vec{u}$ of $u$. This leads us to
define:

\begin{definition}
  \label{def:saturated-polynomial}
  We call the plane polynomial $u = \vec{u} + \eps u_0 \in \DH[t]$
  \emph{saturated} if all real zeros of $\rgcd(\vec{u})$ are of even
  multiplicity. We say that $u$ is \emph{minimally saturated,} if there is no
  polynomial $f \in \R[t]$ of positive degree, such that $u/f$ is a saturated
  polynomial. If $u$ is not saturated, there is a unique monic polynomial $\ell
  \in \R[t]$ of minimal degree, the \emph{saturating factor of $u$,} such that
  $u\ell$ is saturated.
\end{definition}

\begin{example}
  \label{ex:3}
  Consider the polynomial $g = t^2(t-1)^3(t^2+1)$ of Example~\ref{ex:2}. The
  plane polynomial $u = gt\qi + g(t-1)\qj + g(t-2)\qk + \eps t^8$ is not
  saturated as $\rgcd(\vec{u}) = g$ has a real linear factor of odd degree. The
  saturating factor is $\ell = t-1$.
\end{example}

\begin{proposition}
  \label{prop:kinematic-saturated}
  Equation~\eqref{eq:8} has a solution for $P \in \H[t]$ if and only if
  $\vec{u}$ is the vector part of a kinematic and saturated plane
  polynomial~$u$.
\end{proposition}

\begin{proof}
  We already argued that $u$ being kinematic is a necessary condition for
  solutions to exist. The plane polynomial $u$ also needs to be saturated: Set
  $g \coloneqq \rgcd(\vec{u})$ and consider a real polynomial factor $f$ of $g$
  of degree at most two. By the AB-Lemma, cf. \cite[Proposition~2.1]{cheng16} or
  \cite[Lemma~2]{li25}, there are two possibilities:
  \begin{enumerate}
  \item The polynomial $f$ is a factor of $P\qk$ or of $\cj{P}$. But then $f$ is
    a factor of both, $P\qk$ and $\cj{P}$, and appears with quadratic
    multiplicity in~$g$.
  \item There is a polynomial $F \in \H[t]$ such that $F$ is a right factor of
    $P\qk$, $\cj{F}$ is a left factor of $\cj{P}$, and $f = \Norm{F}$. Because
    $\Norm{F}$ is non-negative all zeros of $f$ are complex or of even
    multiplicity.
  \end{enumerate}
  Both cases lead to $u$ being saturated.

  Conversely, if $u$ is kinematic and saturated, Lemma~\ref{lem:2} provides a
  solution polynomial $P \in \H[t]$ with $\ell = 1$.
\end{proof}

\subsection{The Dual Part}
\label{sec:dual}

We continue our study of Equation~\eqref{eq:7} for a given plane polynomial $u =
\vec{u} + \eps u_0$. Denote the yet undetermined motion polynomial by $C = P +
\eps D$. Equation~\ref{eq:7} is equivalent to the following system of equations
for quaternionic polynomials:
\begin{align}
  P \qk \cj{P} &= h\vec{u}, \label{eq:10}\\
  P \qk \cj{D} - D \qk \cj{P} &= hu_0, \label{eq:11}\\
  P \cj{D} + D \cj{P} &= 0. \label{eq:12}
\end{align}
The Equation~\eqref{eq:12} encodes the Study condition and guarantees that $C$
is a motion polynomial. In Section~\ref{sec:primal} we discussed how to solve
Equation~\eqref{eq:10} for the primal part $P$ under the condition $\deg h = 0$.
The extension to the case $\deg h \ge 1$ via Lemma~\ref{lem:1} is
straightforward. Thus, we will mostly focus on Equations \eqref{eq:11} and
\eqref{eq:12} for given $P$, $u_0$, and $h$. It is a system of linear equations
to be solved for $D$ over the ring $\H[t]$. Comparing coefficients of $t$ on
both sides of \eqref{eq:11} and \eqref{eq:12} results in a system of linear
equations for the undetermined real quaternionic coefficients of $D$ which can
readily be solved in concrete examples. It is, however, difficult to make
statements about existence of solutions for given $h$ or minimality of the
solution degree or to understand the structure of solution. Therefore, we do not
convert \eqref{eq:11} and \eqref{eq:12} into a system of linear equations over
$\R$ but instead exploit the rich algebraic environment provided by the
polynomial rings $\H[t]$ and $\R[t]$. Our aim is a characterization of
solubility and of solutions of minimal degree.

We continue with a lemma on existence of solutions. By the end of this section
it will turn out that the solutions described in this lemma are of minimal
degree.

\begin{lemma}
  \label{lem:3}
  Given a reduced kinematic plane polynomial $w$ with minimal saturating factor
  $\ell$, the system of equations~\eqref{eq:10}--\eqref{eq:12} with $u \coloneqq
  \ell w$ and $h \coloneqq \rgcd(\vec{u}/\ell^2) = \rgcd(\vec{w})/\ell$ has a
  solution of degree $\frac{1}{2}(\deg u + \deg h)$.
\end{lemma}

\begin{proof}
  With $\vec{v} \coloneqq \vec{w}/(h \ell)$, Lemma~\ref{lem:2} ensures existence
  of $Q \in \H[t]$ such that $Q \qk \cj{Q} = \vec{v}$. Moreover, $Q$ has no real
  polynomial factor of positive degree.
  The primal part $P$ of $C$ necessarily equals $P = Qh\ell$ and, indeed,
  we have
  \begin{equation*}
    P \qk \cj{P} =(Qh\ell) \qk \cj{(Qh\ell)} =
    h^2 \ell^2 Q \qk \cj{Q} = h^2 \ell^2 \vec{v} = h\ell\vec{w} = h\vec{u}.
  \end{equation*}

  The degree of $P$ equals $\deg P = \frac{1}{2}(\deg u + \deg h)$. Now we need
  to show existence of $D = d_0 + d_1\qi + d_2\qj + d_3\qk \in \H[t]$ with $\deg
  D \le \deg P$ that solves \eqref{eq:11} and \eqref{eq:12}. Substituting $P =
  Qh\ell$ yields
  \begin{align*}
    Q\qk\cj{D} - D\qk\cj{Q} &= w_0,\\
    Q\cj{D} + D\cj{Q} &= 0.
  \end{align*}
  Writing $Q = q_0 + q_1\qi + q_2\qj + q_3\qk$ and $\vec{v} = v_1\qi + v_2\qj +
  v_3\qk$, we have
  \begin{equation}
    \label{eq:13}
    v_1 = 2(q_0q_2 + q_1q_3),\quad
    v_2 = 2(-q_0q_1 + q_2q_3),\quad
    v_3 = q_0^2 - q_1^2 - q_2^2 + q_3^2.
  \end{equation}
  We now assume $q_0^2 + q_3^2 \neq 0$. This is no loss of generality because a
  suitable change of coordinates in the fixed frame (ie., a suitable
  multiplication of $Q$ with a generic quaternion from the left) can ensure
  this. A straightforward computation provides the solution
  \begin{equation}
    \label{eq:14}
    d_0 = \frac{v_2d_1 - v_1d_2 - q_3 w_0}{2(q_0^2+q_3^2)},\quad
    d_3 = -\frac{v_1d_1 + v_2d_2 - q_0 w_0}{2(q_0^2+q_3^2)}
  \end{equation}
  over the field $\R(t)$ of rational functions. All solutions over $\R(t)$ are
  of the shape \eqref{eq:14} but they are not unique because $d_1$ and $d_2$ can
  be chosen freely. Our task is to determine polynomials $d_1$ and $d_2$ such
  that \eqref{eq:14} is polynomial as well. This is the case if and only if
  $d_1$ and $d_2$ solve the system of linear equations
  \begin{equation}
    \label{eq:15}
    \begin{aligned}
      v_2 d_1 - v_1d_2 &= q_3 w_0,\\
      v_1 d_1 + v_2d_2 &= q_0 w_0
    \end{aligned}
  \end{equation}
  over the polynomial ring $R \coloneqq \R[t] / (q_0^2 + q_3^2)$. One solution to
  \eqref{eq:15} is
  \begin{equation}
    \label{eq:16}
    d_1 = -\tfrac{1}{2}w_0 v_3^{-1} q_2,\quad
    d_2 = \tfrac{1}{2}w_0 v_3^{-1} q_1
  \end{equation}
  provided that the inverse of $v_3$ in the ring $R$ exists. Indeed, this can be
  assumed without loss of generality: The inverse $v_3^{-1}$ exists if and only
  if $\gcd(v_3,q_0^2+q_3^2) = 1$. By \eqref{eq:13}, this is the case if and only
  if $\gcd(q_1^2+q_2^2,q_0^2+q_3^2) = 1$. Since $Q$ is reduced this is a
  coordinate-dependent condition and, once more, can be removed by a suitable
  change of coordinates in the fixed frame.

  The solution obtained from \eqref{eq:16} is of degree at most $\deg(q_0^2 +
  q_3^2) - 1 = \deg u - \deg h - 2\deg \ell - 1$. For large degree of $u$ and
  small degrees of $h$ and $\ell$, this violates the condition $\deg D \le \deg
  P = \frac{1}{2}(\deg u + \deg h)$. Thus, some additional work is
  needed.

  Observe that
  \begin{equation*}
    d_1 = q_3,\quad d_2 = -q_0
  \end{equation*}
  is a solution of the homogeneous system corresponding to \eqref{eq:15}.
  Therefore, more solutions to \eqref{eq:15} can be written as
  \begin{equation*}
    d_1 = -\tfrac{1}{2}w_0 v_3^{-1} q_2 + \lambda q_3,\quad
    d_2 = +\tfrac{1}{2}w_0 v_3^{-1} q_1 - \lambda q_0
  \end{equation*}
  for arbitrary $\lambda \in R$. By Lemma~\ref{lem:4} below, with $p = q_0^2 +
  q_3^2$, $x = d_1$, $y = d_2$, $a = q_3$ and $b = -q_0$, there exists $\lambda
  \in R$ such that the desired degree bound can be achieved. We still need to
  argue that the conditions to apply Lemma~\ref{lem:4} are met. These are
  \begin{equation}
    \label{eq:17}
    q_3 \neq -q_0,\quad
    \gcd(q_3,q_0^2+q_3^2) = 1,\quad
    \gcd(q_0,q_0^2+q_3^2) = 1.
  \end{equation}
  Since $Q$ is reduced, we can repeat an already used argument and say that a
  proper choice of coordinates in the fixed frame will guarantee \eqref{eq:17}.
\end{proof}

\begin{lemma}
  \label{lem:4}
  Given $n \in \mathbb{N}$ and $p \in \R[t]$ with $\deg p = 2n$, denote by $R =
  \R[t]/(p)$ the ring of polynomials in $\R[t]$ modulo $p$. For any $a$, $b$,
  $x$, $y \in R$ with $\deg a \le n$, $\deg b \le n$, $a \neq b$, $\gcd(a,p) =
  \gcd(b,p) = 1$ there exists $\lambda \in R$ such that
  \begin{equation*}
    \deg(x + \lambda a) \le n
    \quad\text{and}\quad
    \deg(y + \lambda b) \le n.
  \end{equation*}
\end{lemma}

\begin{proof}
  We identify $R$ with the vector space $\R^{2n}$. For any $z$, $c \in R$,
  the map $F_{z,c}\colon \R^{2n} \to \R^{2n}$, $\lambda \mapsto z + \lambda c$
  is affine. If $\gcd(p,c) = 1$, we have
  \begin{equation*}
    \kernel(F_{0,c}) = \{\lambda \in \R^{2n} \mid \lambda c = 0\} = \{ 0 \}.
  \end{equation*}
  Thus, $F_{0,a}$, $F_{0,b}$, $F_{x,a}$, and $F_{y,b}$ are all bijections.

  Denote by $U \subset \R^{2n}$ the subspace of polynomials of degree at most
  $n$. It is of dimension $n + 1$. We have to show that $F_{x,a}^{-1}(U) \cap
  F_{y,b}^{-1}(U)$ is not empty.

  Consider the linear subspaces $F_{0,a}^{-1}(U)$ and $F_{0,b}^{-1}(U)$. Both
  are of dimension $n+1$ and we claim that they are not identical. Indeed, the
  linear maps $U \to F_{0,a}^{-1}(U)$, and $U \to F_{0,b}^{-1}(U)$ are given by
  $u \mapsto ua^{-1}$, $u \mapsto ub^{-1}$ where $a^{-1}$, $b^{-1}$ are the
  respective inverse elements in $R$. Thus, $F_{0,a}^{-1}(U)$
  and $F_{0,b}^{-1}(U)$ are identical if and only if for all $u \in U$ we have
  $ua^{-1} = ub^{-1}$. Since $U$ contains invertible elements, this implies
  $a^{-1} = b^{-1}$ and consequently $a = b$. This is excluded by assumption and
  $F_{x,a}^{-1}(U) \cap F_{yb}^{-1}(U)$ is not empty. Any polynomial $\lambda$
  from this intersection solves our problem.
\end{proof}

\begin{remark}
  \label{rem:uniqueness}
  The construction in the proof of Lemma~\ref{lem:4} leaves two degrees of
  freedom in the choice of $\lambda$ as two vector subspaces of dimension $n+1$
  in $\R^{2n}$ intersect in a subspace of dimension at least two. In fact, they
  intersect precisely in a subspace of this dimension since $ua^{-1} = ub^{-1}$
  implies $u = 0$ unless $a^{-1} = b^{-1}$. Thus, there are also two free
  parameters in the choice of the polynomials $d_1$ and $d_2$ in the proof of
  Lemma~\ref{lem:3}. This can also be seen as follows: Given a solution $C = P +
  \eps D$, we can right-multiply it with a dual quaternion of the shape $e
  \coloneqq 1 + \eps (e_5\qi + e_6\qj) \in \DH$. Because of $\ej{e} \qk \cj{e} =
  \qk$, $C$ and $Ce$ are both solutions of the same degree.

  Note however, that the construction of Lemma~\ref{lem:4} produces polynomials
  $d_1$, $d_2$ such that $\deg D \le \deg Q$ while Equation~\eqref{eq:14} only
  demands the degree bound $\deg\vec{v} + \deg D \le \deg u + \deg Q$. This
  allows to replace $U$ in the proof of Lemma~\ref{lem:4} by the subspace of
  polynomials of degree at most $n + 2\deg \ell h$ and accounts for a total of $2
  + 2\deg \ell h$ free parameters in the selection of $\lambda$. This observation
  will be important when we discuss uniqueness of solutions.
\end{remark}

\begin{example}
  \label{ex:4}
  We consider the plane polynomial $u = \vec{u} + \eps u_0 = u_1\qi + u_2\qj +
  u_3\qk + \eps u_0 = h(v_1\qi + v_2\qj + v_3\qk) + \eps
  u_0$ where
  \begin{gather*}
    u_0 = 2t^2 - 14t + 20,\quad
    v_1 = 2(t^2-3t+1),\quad
    v_2 = -2(t^2-2t+2),\\
    v_3= -(t^2+2t-4),\quad\text{and}\quad
    h = t^2-6t+10.
  \end{gather*}
  The plane polynomial is kinematic, reduced, and saturated, hence $\ell = 1$.
  The vector part of $u$ has the non-trivial $\rgcd(\vec{u}) = h$. Using the
  method of \cite[Lemma~2.3]{schroecker24} we solve
  \begin{equation*}
    Q\qk\cj{Q} = v_1\qi + v_2\qi + v_3\qk
  \end{equation*}
  for $Q$ and obtain
  \begin{equation*}
    Q = 1 + t\qi + (1-t)\qj + (t-2)\qk.
  \end{equation*}
  This gives us the primal part
  \begin{equation*}
    P \coloneqq Qh =
    (t^2-6t+10)(1 + t\qi + (1-t)\qj + (t-2)\qk).
  \end{equation*}
  Indeed, we have $P\qk\cj{P} = h^2 Q\qk\cj{Q} = h\vec{u}$.

  In order to find the dual part $D$ of the sought motion polynomial $C = P +
  \eps D$, we compute the inverse of $v_3$ in the ring $R = \R[t]/(q_0^2+q_3^2)$
  by the extended Euclidean algorithm. It is $v_3^{-1} =
  \frac{2}{15}t-\frac{1}{3}$. Now we use \eqref{eq:16} to compute one solution
  for $d_1$ and $d_2$ in $R$:
  \begin{equation*}
    d_1 = \tfrac{8}{15}t - \tfrac{2}{3},\quad
    d_2 = \tfrac{3}{5}t - \tfrac{1}{3}.
  \end{equation*}
  Since the degree of $d_1$ and $d_2$ is already low enough, we do not need to
  reduce it further by adding suitable multiples of $q_3$ and $-q_0$,
  respectively (Lemma~\ref{lem:4}). Using Equation~\eqref{eq:14} we find
  \begin{equation*}
    d_0 = -\tfrac{32}{15}t + \tfrac{13}{3},\quad
    d_3 = \tfrac{1}{15}t + 2.
  \end{equation*}
  This gives the solution polynomial
  \begin{multline}
    \label{eq:18}
    C =
    (\qi-\qj+\qk)t^3
    +(-6\qi+7\qj-8\qk+1)t^2\\
    -(6 - 10\qi+16\qj-22\qk+\tfrac{1}{15}\eps(32-8\qi-9\qj-1\qk))t\\
    +10+10\qj-20\qk+\tfrac{1}{3}\eps(13-2\qi-1\qj+6\qk).
  \end{multline}
  Further solutions can be constructed from \eqref{eq:18} by multiplying $C = P
  + \eps D$ from the right with a polynomial of the shape $E = 1 + \eps (e_5 \qi
  + e_6 \qj)$. If $E$ is of degree at most two, this will not increase the
  degree of $C$. As predicted by Remark~\ref{rem:uniqueness}, there are six
  degrees of freedom, the real coefficients $e_5$ and~$e_6$.
\end{example}

The motion $[C]$ with $C$ given by \eqref{eq:18} will be visualized later in
Figure~\ref{fig:1}.

\begin{lemma}
  \label{lem:5}
  Given a saturated kinematic plane polynomial $u \in \DH[t]$, the system of
  equations~\eqref{eq:10}--\eqref{eq:12} has no solution if $h$ is a proper
  factor of $\rgcd(\vec{u})$.
\end{lemma}

\begin{proof}
  Since otherwise no solutions exist, we can assume that all real linear factors
  of $h$ are of even multiplicity. We set $g \coloneqq \rgcd(\vec{u})$ and $k
  \coloneqq g/h$. By assumption $k$ is a polynomial of positive degree and all
  its real linear factors are of even multiplicity. By Lemma~\ref{lem:1}, there
  exist polynomials $H$, $K \in \C[t]$ such that $h = H\cj{H}$ and $k =
  K\cj{K}$. We also set $G \coloneqq KH$.

  In order to find a motion polynomial solution to Equations
  \eqref{eq:10}--\eqref{eq:12}, we follow the steps of Lemma~\ref{lem:3}. The
  primal part is necessarily of the shape $P = QG\cj{H} = QKH\cj{H} = QKh$ with
  $Q \in \H[t]$ and $\gcd(Q\cj{Q}, \rgcd(\vec{u})) = 1$. In particular, $Q$ is
  free of real polynomial factors of positive degree. Plugging this into
  Equations~\eqref{eq:11} and \eqref{eq:12} yields
  \begin{align*}
    (QK)\qk\cj{D} - D\qk\cj{(QK)} - u_0 &= 0,\\
    (QK)\cj{D} + D\cj{(QK)} &= 0.
  \end{align*}
  With $\vec{v} \coloneqq (QK)\qk\cj{(QK)} = kQ\qk\cj{Q} = k(v_1\qi + v_2\qj +
  v_3\qk)$ and $Q = q_0 + q_1\qi + q_2\qj + q_3\qk$ the system of equations
  \eqref{eq:15} becomes
  \begin{equation}
    \label{eq:19}
    \begin{aligned}
      k(v_2d_1 - v_1d_2) &= q_3u_0, \\
      k(v_1d_1 + v_2d_2) &= q_0u_0.
    \end{aligned}
  \end{equation}
  It is to be solved over $R = \R[t]/(q_0^2+q_3^2)$. Since $Q$ has no polynomial
  factors of positive degree, we can assume without loss of generality that
  $\gcd(q_0,q_0^2+q_3^2) = \gcd(q_3,q_0^2+q_3^2) = 1$ by left-multiplication
  with a generic quaternion. If this quaternion is selected to be of unit norm,
  then $u_0$ will not change and we can also assume $\gcd(u_0,q_0^2+q_3^2) = 1$.
  Hence, the right-hand sides of \eqref{eq:19} are invertible in $R$.

  By definition of $Q$ and $v$, the relations \eqref{eq:13} are to be replaced by
  \begin{equation*}
    kv_1 = 2(q_0q_2+q_1q_3),\quad
    kv_2 = 2(-q_0q_1+q_2q_3),\quad
    kv_3 = q_0^2 - q_1^2 - q_2^2 + q_3^2.
  \end{equation*}
  But then
  \begin{multline*}
    q_2(q_0^2 + q_3^2) =
    q_2q_0^2+q_0q_1q_3 - q_0q_1q_3+q_2q_3^2 \\=
    (q_0q_2+q_1q_3)q_0 + (-q_0q_1+q_2q_3)q_3 =
    \tfrac{1}{2}k(v_1 q_0 + v_2 q_3).
  \end{multline*}
  Assuming once more without loss of generality that $\gcd(k,q_2) = 1$ we infer that
  $k$ divides $q_0^2+q_3^2$ as well as the left-hand sides of \eqref{eq:19}.
  Since $\deg k \ge 1$, this means that the left-hand sides are not invertible
  in $R$. Thus, there exists no solution of \eqref{eq:19} and, consequently, no
  solution to equations~\eqref{eq:10}--\eqref{eq:12} either.
\end{proof}

\begin{lemma}
  \label{lem:6}
  If the system of equations~\eqref{eq:10}--\eqref{eq:12} has a solution and $f$
  is a quadratic $\R$-irreducible factor of $h/\gcd(\rgcd(\vec{u}),h)$ whose
  multiplicity as factor of $h$ is one, then there exists a solution to
  \eqref{eq:10}--\eqref{eq:12} with $h$ replaced by $h/f$.
\end{lemma}

\begin{proof}
  We denote the solution to \eqref{eq:10}--\eqref{eq:12} by $P$ and $D$ and set
  $C \coloneqq P + \eps D$. The idea of our proof is to show existence a motion
  polynomial $E$ that is a right factor of $C$, ie. $C = \tilde{C}E$, and
  satisfies $\ej{E}\qk\cj{E} = f\qk$. Then, because of
  \begin{equation*}
    \ej{\tilde{C}} \qk \cj{\tilde{C}} =
    \tfrac{1}{f}\ej{(\tilde{C}E)} \qk \cj{(\tilde{C}E)} =
    \tfrac{1}{f}\ej{C}\qk\cj{C} =
    \tfrac{h}{f} u,
  \end{equation*}
  the primal and dual parts of $\tilde{C}$ provide the claimed solution.

  There exist $F = f_0 + f_3\qk \in \C[t] \setminus \R[t]$ and $Q \in \H[t]$
  such that $P = QF$ and $F\cj{F} = f$. Moreover, $f$ is not a factor of $Q$ by
  assumption. Therefore, \cite[Lemma~2 and Lemma~3]{hegedus13} implies existence
  of a linear motion polynomial $E = F + \eps K$ that is a right factor of $C =
  P + \eps D$, ie. $C = \tilde{C}E$ for some motion polynomial $\tilde{C}$, such
  that $E\cj{E} = f$. This later condition implies that $E$ is of the shape $E =
  f_0 + f_3\qk + \eps(e_5\qi + e_6\qj)$ with real polynomials $e_5$, $e_6$.
  Hence, $\ej{E}\qk\cj{E} = (f_0^2+f_3^2)\qk = f\qk$. Finally, the norm polynomials
  of $C$ and $E$ are both real. Because the norm is multiplicative, the same is
  true for $\Norm{\tilde{C}}$ and $\tilde{C}$ is a motion polynomial.
\end{proof}

\begin{example}
  \label{ex:5}
  Lets study again the plane polynomial $u$ of Example~\ref{eq:4}. The motion
  polynomial
  \begin{multline}
    \label{eq:20}
    \hat{C} = (\qi-\qj+\qk-\eps(2-\qi-\qj))t^4
    +(2-5\qi+8\qj-8\qk+\eps(13-8\qi-8\qj+\qk))t^3 \\
    -(14-3\qi+22\qj-21\qk+\tfrac{1}{15}\eps(437-383\qi-324\qj+104\qk))t^2 \\
    +(32+16\qi+20\qj-14\qk+\tfrac{1}{15}\eps(306-589\qi-207\qj+302\qk))t \\
    -20-10\qi-10\qk - \tfrac{1}{3}\eps(24-91\qi+32\qj+43\qk)
  \end{multline}
  satisfies the equality $\ej{\hat{C}} \qk \cj{\hat{C}} = hu$ where $h =
  (t^2+1)g$ and $g = \rgcd(\vec{u})$. Thus, Lemma~\ref{lem:6} suggests that
  there is a linear polynomial $E = t - \qk + \eps(e_5\qi + e_6\qj)$ or $E = t +
  \qk + \eps(e_5\qi + e_6\qj)$ in $\DH[t]$ that is a right factor of $C$.
  Indeed, with
  \begin{equation*}
    E = t - \qk + \eps ((\qi - \qj)t + \qi + \qj)
  \end{equation*}
  we have $\hat{C} = CE$ where $C$ is the motion polynomial of
  Equation~\eqref{eq:18}.
\end{example}

The motion $[\hat{C}]$ with $\hat{C}$ given by \eqref{eq:20} will be visualized
later in Figure~\ref{fig:1}.

\begin{lemma}
  \label{lem:7}
  If the system of equations~\eqref{eq:10}--\eqref{eq:12} has a solution and
  there exists a linear or quadratic $\R$-irreducible real polynomial $f$ such
  that $f^2 \in \R[t]$ is a factor of $h/\gcd(\rgcd(\vec{u}),h)$, then there
  exists $m \in \mathbb{N}$, $m \ge 1$ and a solution to
  \eqref{eq:10}--\eqref{eq:12} with $h$ replaced by $h/f^{2m}$.
\end{lemma}

\begin{proof}
  Denote the solution to \eqref{eq:10}--\eqref{eq:12} by $C = P + \eps D$. The
  primal part is of the shape $P = Qf^m$ for some positive integer $m$ and $Q
  \in \H[t]$ with $\gcd(f,\rgcd(Q))=1$.

  We claim existence of a motion polynomial $E = f^m + \eps F \in \DH[t]$ that
  is a right factor of $C$ and satisfies
  \begin{equation}
    \label{eq:21}
    (f^m - \eps F)\qk(f^m + \eps \cj{F}) = f^{2m}\qk.
  \end{equation}
  As already argued at the beginning of the proof of Lemma~\ref{lem:6}, this
  would imply the lemma's statement.

  Using polynomial division, we find $D'$, $R \in \DH[t]$ such that $D = D'f^m +
  R$. Lets set $q \coloneqq Q\cj{Q}$ and $\tilde{F} \coloneqq \cj{Q}Rq^{-1}$
  where $q^{-1}$ is the inverse of $q$ in the polynomial ring $\R[t]/(f^m)$. It
  exists because of $\gcd(f,\rgcd(Q)) = 1$. In addition, we define $F$ to be
  the unique representative of $\tilde{F}$ in the ring of quaternions with
  coefficients in $\R[t]/(f^m)$ that is of degree less than $f^m$ ($F$ is the
  remainder of $\tilde{F}$ divided by $f^m$). This implies that $QF = Lf^m +
  R$ for some $L \in \H[t]$. With $K \coloneqq D' - L$ we now have
  \begin{multline*}
    C = Qf^m + \eps D
    = Qf^m + \eps (D'f^m - Lf^m + Lf^m + R)\\
    = Qf^m + \eps(Kf^m + QF)
    = (Q + \eps K)(f^m + \eps F).
  \end{multline*}
  The norm of $C$ equals $Q\cj{Q}f^{2m}$ whence the norm of $f^m + \eps F$
  equals $f^{2m}$. Therefore, $f^{2m} + \eps F$ is a motion polynomial, that is,
  $F = f_{21}\qi + f_{22}\qj + f_{23}\qk$ with $f_{21}$, $f_{22}$, $f_{23} \in
  \R[t]$.

  We still need to show that $f^m + \eps F$ satisfies \eqref{eq:21}. This is the
  case if and only if $f_{23} = 0$. Recall that $f^{2m}$ is a factor of
  $\ej{C}\qk\cj{C}$. We compute
  \begin{equation*}
    \begin{aligned}
      \ej{C}\qk\cj{C}
      &= (Q - \eps K)(f^m - \eps F)\qk(f^m + \eps \cj{F})(\cj{Q} + \eps\cj{K}) \\
      &= f^{2m}(Q\qk\cj{Q} + \eps(Q\qk\cj{K}-K\qk\cj{Q})) + \eps f^m Q(\qk\cj{F} - F\qk)\cj{Q} \\
      &= f^{2m}(Q\qk\cj{Q} + \eps(Q\qk\cj{K}-K\qk\cj{Q})) + 2\eps f^m f_{23} q
    \end{aligned}
  \end{equation*}
  and infer that $f^m$ is a factor of $f_{23}q$. Since $f$ is not a factor of
  $q$, $f^m$ is a factor of $f_{23}$. But since $\deg f_{23} \le \deg F <
  \deg f^m$, we have $f_{23} = 0$.
\end{proof}

\begin{example}
  \label{ex:6}
  We once more look at the plane polynomial $u$ of Example~\ref{ex:4}. The
  motion polynomial
  \begin{multline*}
    \tilde{C} = (1-6\qi+7\qj-8\qk-(2-\qi-\qj)\eps)t^4 \\
    -(6-10\qi+16\qj-22\qk-\tfrac{1}{15}(163-112\qi-111\qj+16\qk)\eps)t^3 \\
    +(10+10\qj-20\qk+\tfrac{1}{3}(-68+73\qi+62\qj-15\qk)\eps)t^2 \\
    +(16-38\qi-14\qj+16\qk)\eps t - (10-30\qi+10\qj+10\qk)\eps
  \end{multline*}
  satisfies the equality $\ej{\tilde{C}} \qk \cj{\tilde{C}} = hu$ where $h =
  t^4g$ and $g = \rgcd(\vec{u})$. Thus, Lemma~\ref{lem:7} suggests that with $f
  = t^2$ there is a quadratic motion polynomial $E = f + \eps(e_5\qi + e_6\qj)
  \in \DH[t]$ that is a right factor of $\tilde{C}$. Indeed, with
  \begin{equation*}
    E = t^2 + ((\qi - \qj) t + \qi + \qj)\eps
  \end{equation*}
  we have $\tilde{C} = CE$ where $C$ is the motion polynomial of~\eqref{eq:18}.
\end{example}

We now summarize the findings of this section in a proposition:

\begin{proposition}
  \label{prop:plane-polynomial}
  Given a kinematic, minimally saturated plane polynomial $u$ with minimal
  saturating factor $\ell$ and $h \in \R[t]$,
  equations~\eqref{eq:10}--\eqref{eq:12} have a solution if and only if $h$ is a
  multiple of $\rgcd(\ell\vec{u})/\ell^2$ and $h\ell u$ is still saturated.
  Solutions of minimal degree are obtained for $h = \rgcd(\ell\vec{u})/\ell^2$.
\end{proposition}

\begin{proof}
  The proof is a combination of Lemmas~\ref{lem:3}, \ref{lem:5}, \ref{lem:6},
  and \ref{lem:7}.

  If $u$ is kinematic and minimally saturated and $h = \rgcd(\vec{u})$, a
  solution $C = P + \eps D$ exists by Lemma~\ref{lem:3}. If $h$ is a multiple of
  $\rgcd(\vec{u})$ such that $hu$ is still saturated, further solutions are found
  by right multiplication of $C$ with a suitable motion polynomial $E \in
  \DH[t]$ such that $\Norm{E} = h / \rgcd(\vec u)$.

  Conversely, if $C = P + \eps D$ is a solution for some $h \in \R[t]$, then
  $hu$ is necessarily kinematic and saturated. Using Lemma~\ref{lem:6} and
  \ref{lem:7} we can construct further solutions by splitting off irreducible
  factors of multiplicity one or linear or quadratic factors of even multiplicity
  from $h$ until, by Lemma~\ref{lem:5}, we find a solution with $h =
  \rgcd(\vec{u})$. Thus, $h$ is a multiple of $\rgcd(\vec{u})$.

  Lemmas~\ref{lem:3} and \ref{lem:5} together with Lemma~\ref{lem:2} also imply
  the statement about the minimal degree.
\end{proof}

\section{Motions of Minimal Degree}
\label{sec:minimal}

In Section~\ref{sec:system} we have discussed the equation $\ej{C} \qk \cj{C} =
hu$ for a kinematic plane polynomial $u \in \DH[t]$, a real polynomial $h \in
\R[t]$, and a motion polynomial $C \in \DH[t]$. We formulated statements on
solubility and on solutions of minimal degree in terms of real and quaternionic
polynomial algebra. The purpose of this section is to translate our findings
into the language of kinematics and to state them in one central theorem. Let us
start by transferring the concept of a kinematic plane polynomial.

\begin{definition}
  \label{def:kinematic-torse}
  A rational torse $[u]$ is called \emph{kinematic} if it can be represented by
  a kinematic plane polynomial~$u$.
\end{definition}

\begin{remark*}
  Only kinematic torses can arise as trajectories of a plane under a rational
  motion. They can be characterized among all rational torses as those having a
  rational Gauss map $n \coloneqq \vec{u}/\sqrt{\Norm{u}}$. Thus, the tangent
  planes of a cone of revolution form the planes of a kinematic rational torse
  while the tangent planes of a general quadratic cone do not. The former
  obviously arise trajectory of a rational motion, the rotation around the
  cone's axis. For the latter, Equation~\eqref{eq:7} and also
  Equation~\eqref{eq:8} have no solution.
\end{remark*}

The following theorem contains the main results of this article. It uses the
notion of “essential uniqueness” of the rational motion $[C]$. By this we mean
that $[C]$ is unique up to right multiplication with a constant displacement
$[e]$ that fixes the plane $[\qk]$, that is $[\ej{e}\qk\cj{e}] = [\qk]$.

\begin{theorem}
  \label{th:1}
  Given a rational kinematic torse $[u] = [\vec{u} + \eps u_0]$ with Gauss map
  $n = \vec{u}/\sqrt{\Norm{\vec{u}}}$, the following hold true:
  \begin{enumerate}
  \item There exists a rational motion $[C]$ of degree $\deg u - \frac{1}{2}\deg
    n$ with trajectory $[u]$.
  \item The degree of this rational motion $[C]$ is minimal.
  \item The rational motion of minimal degree is essentially unique if and only
    if $\deg u = \deg n$.
  \item All rational motions with trajectory $[u]$ are obtained by composing
    $[C]$ from the right with any rational motion $[E]$ that fixes the plane
    $[\qk]$.\footnote{This implies that $E$ is given by a polynomial of the
      shape $E = e_0 + e_3\qk + \eps (e_5\qi + e_6\qj)$ with $e_0$, $e_3$,
      $e_5$, $e_6 \in \R[t]$.}
  \end{enumerate}
\end{theorem}

\begin{proof}
  We assume that the rational kinematic torse $[u]$ is given by a reduced
  kinematic plane polynomial $u$. Denote by $\ell$ the minimal saturating factor
  of $u$.

  Now, existence of the rational motion $[C]$ and minimality of its degree is a
  consequence of Proposition~\ref{prop:plane-polynomial}. With $h =
  \rgcd(\ell\vec{u})/\ell$ the Gauss map $n$ is of degree $\deg n = \deg u -
  \deg h - \deg \ell$ and the minimal motion degree is
  \begin{equation*}
    \deg C = \tfrac{1}{2}(\deg u + \deg \ell + \deg h) = \tfrac{1}{2}(2\deg u -
    \deg n) = \deg u - \tfrac{1}{2}n.
  \end{equation*}

  The statement on its essential uniqueness follows from
  Remark~\ref{rem:uniqueness}.

  Finally, any solution $[\hat{C}]$ can be reduced to the minimal degree
  solution $[C]$ by splitting of suitable motion polynomial right factors $E$ as
  shown in Lemmas~\ref{lem:6} and \ref{lem:7}. Thus, $[\hat{C}]$ is conversely
  obtained from $[C]$ by composition with $[E]$ from the right.
\end{proof}

\begin{figure}
  \centering
  \includegraphics[]{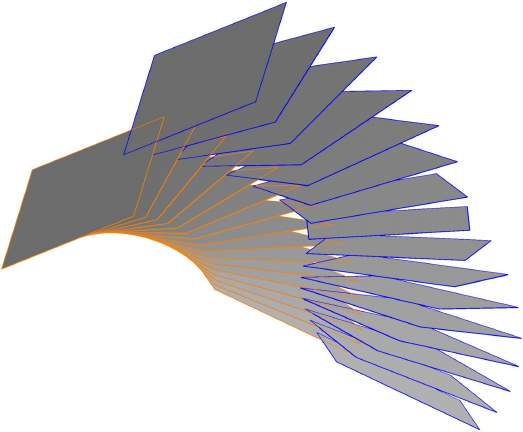}
  \caption{Two different minimal degree motions with identical plane
    trajectory.}
  \label{fig:1}
\end{figure}

In order to actually compute rational motions of minimal degree we may proceed
as in Examples~\ref{ex:4} and \ref{ex:5}. Figure~\ref{fig:1} visualizes the
trajectories of the plane $[\qk]$ for the rational motion $[C]$ given by
\eqref{eq:18} in Example~\ref{ex:4} and the rational motion $[\hat{C}]$ given by
\eqref{eq:20} in Example~\ref{ex:5}. The positions of the moving plane are
represented by moving rectangles in both cases. The overlap of the two
rectangles confirms that they lie in the same plane but the overlapping region
changes size so that the two motions are indeed essentially different.

\begin{figure}
  \centering
  \includegraphics[]{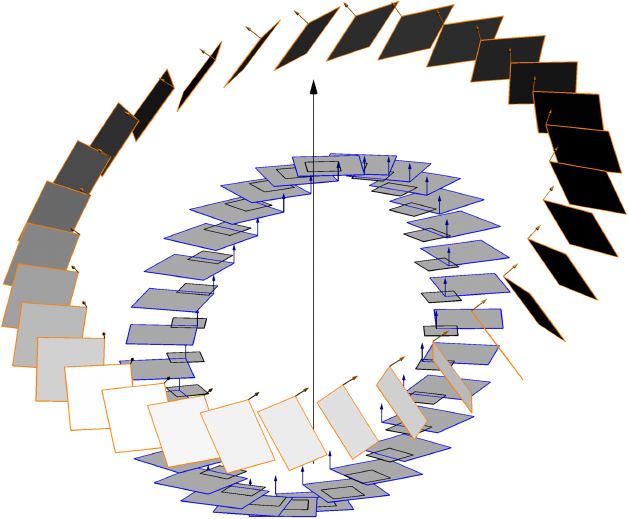}
  \caption{Plane trajectories of the Darboux motion.}
  \label{fig:2}
\end{figure}

\begin{example}
  As a concluding example, lets consider the Darboux motion, the famous rational
  space motion with only planar trajectories \cite[Ch.~9, \S~3]{Bottema90}. It
  can be thought of as the composition of a Cardan motion in a horizontal plane
  with a harmonic oscillation along a vertical axis. Figure~\ref{fig:2}
  visualizes the motion of three planes undergoing a Darboux motion. The first
  trajectory is that of a generic plane where both, the torse and the Gauss map
  are of degree two. Hence, the Darboux motion is minimal for this torse. The
  second trajectory is that of a horizontal plane. We clearly observe that the
  planes' normals have a fixed direction whence the degree of the Gauss map is
  zero. The Darboux motion is minimal but not unique. It can be composed from
  the right with a suitable translational motion with zero vertical component.
  This is illustrated by the smaller frames. Clearly, they lie in the same plane
  but they also change their relative position with respect to the larger frame.
\end{example}

\section{Conclusion and Future Research}
\label{sec:conclusion}

We have presented a complete discussion of the problem to create a prescribed
rational torse $[u]$ as plane trajectory of a rational motion. The torse
necessarily needs to be kinematic, that is, have a possess a rational Gauss map.
In general, if $\rgcd(\vec{u}) = 1$ the resulting motion of minimal degree is
unique up to inessential coordinate transformations. Our constructive proofs
allow to compute solutions of minimal degree and make effective use of the
algebra of real and quaternionic polynomials. All solutions to the problem at
hand are obtained by composing from the right a solution of minimal degree with
rational planar motions that fix the moving plane.

A natural next question is to study the problem of creating a given rational
ruled surface as trajectory of a straight line undergoing a rational motion. As
demonstrated in the pre-print \cite{derin25}, fundamental statements are quite
similar to the torse case but proof details differ.

The algebra of dual quaternions only allows to model the group $\SE$ of rigid
body displacements of Euclidean three-space. It might be a worthy undertaking to
study similar problems in other algebras and kinematic spaces, for example the
creation of one-parametric sets of circles by low degree motions in conformal
kinematics. Possibly, this could unveil relations to rationally
parametrized channel surfaces \cite{Peternell97,Landsmann01}.

We do believe that a thorough understanding of rational motions to create plane
and line trajectories, its restrictions and its degrees of freedom, will
ultimately also pave the way towards further application in engineering sciences.

\section*{Acknowledgment}

Zülal Derin Yaqub was supported by the Austrian Science Fund (FWF) P~33397-N
(Rotor Polynomials: Algebra and Geometry of Conformal Motions) and also
gratefully acknowledges the support provided by the Scientific and Technological Research Council of Türkiye, TUBITAK-2219 -- International Postdoctoral Research Fellowship Program for Turkish Citizens.

\bibliographystyle{elsarticle-num}

\end{document}